\documentclass[10pt]{article}

\usepackage{theorem,amssymb,amsmath}

\usepackage{graphicx}

\usepackage{color}                    
\definecolor{red}{rgb}{1,0,.2}        
\definecolor{cjp}{rgb}{.1,.7,.2}        
\definecolor{fmdc}{rgb}{1,0,.8}        

\newcommand{\A}{{\rm A}}
\newcommand{\B}{{\rm B}}
\newcommand{\C}{{\rm C}}
\newcommand{\D}{{\rm D}}

\advance \topmargin by -\headheight
\advance \topmargin by -\headsep
\evensidemargin \oddsidemargin
\author{
F. M. Dekking \\
Delft University of Technology \\
Faculty EEMCS, P.O.~Box 5031\\
2600 GA Delft, The Netherlands\\
{\tt F.M.Dekking@math.tudelft.nl}
}

\title{\bf Base phi representations and\\ golden mean beta-expansions }

\date{\today}

\def \proof{\noindent{\it Proof:\ \ }}

\def \endpf{{\ \ $\Box$ \medbreak}}

\newtheorem{theorem}{Theorem}[section]

\newtheorem{conjecture}{Conjecture}
\newtheorem{lemma}[theorem]{Lemma}

\theorembodyfont{\rm}

\newtheorem{remark}[theorem]{Remark}

\begin{document}

\maketitle

\begin{abstract}
 In the base phi representation any natural number is written uniquely as a sum powers of the golden mean with digits 0 and 1, where one requires that the product of two consecutive digits is always 0. In this paper we give precise expressions for the those natural numbers for which the $k$th digit is 1, proving two conjectures for $k=0,1$. The expressions are all in terms of generalized Beatty sequences.
\end{abstract}

\bigskip

\section{Introduction}

  Base phi representations were introduced by George Bergman in 1957 (\cite{Bergman}). Base phi representations are also known  as beta-expansions of the natural numbers, with $\beta=(1+\sqrt{5})/2=:\varphi$, the golden mean.\\
  A natural number $N$ is written in base phi if $N$ has the form
  $$N= \sum_{i=-\infty}^{\infty} d_i \varphi^i,\vspace*{-.0cm}$$
  with digits $d_i=0$ or 1, and where $d_1d_{i+1} = 11$ is not allowed.
  Similarly to base 10 numbers, we  write these representations as
  $$\beta(N) = d_{L}d_{L-1}\dots d_1d_0\cdot d_{-1}d_{-2} \dots d_{R-1}d_R.$$

  The base phi representation of a number $N$ is unique (\cite{Bergman}).
  Our main concern will be the distribution of the digit $d_0=d_0(N)$ over the natural numbers $N\in \mathbb{N}$. Several authors have interpreted this in the frequency sense.
  The following result was conjectured by Bergman, and proved in \cite{Hart99}.

\begin{theorem} The frequency of $1$'s in $(d_0(N))$ exists, and
$ \lim_{N\rightarrow\infty} \frac1N \sum_{M=1}^{N} d_0(M)  = \frac1{\varphi+2}=\frac{5-\sqrt{5}}{10}.$
\end{theorem}

A more detailed description, obviously implying the previous theorem, was conjectured  by Baruchel  in 2018 (see A214971 in \cite{oeis}):

\begin{conjecture} Digit $d_0(N)=1$ if and only if $N=\lfloor n\varphi\rfloor +2n +1$ for some natural number $n$, or $N=1$.
\end{conjecture}


Here $\lfloor \cdot \rfloor$ denotes the floor function, and $(\lfloor n\varphi\rfloor)$ is the well known lower Wythoff sequence.
The corresponding result for digit $d_1$ was conjectured  by Kimberling in 2012 (see A054770 in \cite{oeis}):

\begin{conjecture} Digit $d_1(N)=1$ if and only if $N=\lfloor n\varphi\rfloor +2n-1$ for some natural number $n$.
\end{conjecture}

\medskip


 Both conjectures will be proved in Section \ref{sec:main}.  In Section \ref{sec:Beatty}, \ref{sec:morph} and \ref{sec:luc} we introduce some objects and tools used in the proof. Finally Section \ref{sec:general} gives the result for any digit $d_k(N)$ with $k\ge 1$ of the base phi expansion.

 \medskip

 In future work we plan to extend our results to the metallic means, or more generally to arbitrary
 quadratic bases, as defined and analyzed in \cite{Burger}.

\section{Generalized Beatty sequences}\label{sec:Beatty}

  The sequences occurring in the conjectures are sequences $V$ of the type  $V(n) = p(\lfloor n \alpha \rfloor) + q n +r $, $n\ge 1$,  where $\alpha$ is a real number, and $p,q,$ and $r$ are integers. As in \cite{GBS}, we call them \emph{generalized Beatty sequences}.
 \noindent If $S$ is a sequence, we denote its  sequence of first order differences as $\Delta S$, i.e., $\Delta S$ is defined by
$$\Delta S(n) = S(n+1)-S(n), \quad {\rm for\;} n=1,2\dots$$

It is well known (\cite{lothaire}) that the sequence $\Delta(\lfloor n\varphi\rfloor)$  is equal to the Fibonacci word $x_{1,2} = 1211212112\dots$  on the alphabet $\{1,2\}$. More generally, we have  the following simple lemma.

\smallskip

\begin{lemma}\label{lem:diff}{ \rm \bf(\cite{GBS})} Let $V = (V(n))_{n \geq 1}$ be the generalized Beatty
sequence defined by $V(n) = p\lfloor n \varphi \rfloor + q n +r$, and let $\Delta V$ be the
sequence of its first differences. Then $\Delta V$ is the Fibonacci word on the alphabet
$\{2p+q, p+q\}$. Conversely, if $x_{a,b}$ is the Fibonacci word on the alphabet
$\{a,b\}$,  then any $V$ with $\Delta V= x_{a,b}$ is a generalized Beatty sequence
$V=((a-b) \lfloor n \varphi \rfloor)+(2b-a)n+r)$ for some integer $r$.
\end{lemma}

\section{Morphisms}\label{sec:morph}
A morphism is a map from the set of infinite words over an alphabet to itself, respecting the concatenation operation.
The canonical example is the Fibonacci morphism $\sigma$ on the alphabet $\{0,1\}$ given by
$$\sigma(0)= 01, \quad \sigma(1)= 0.$$
A central role in this paper is played by the morphism $\gamma$ on the alphabet  $\{\A,\B,\C,\D\}$ given by
$$\gamma(\A)=\A\B, \quad \gamma(\B)=\C, \quad \gamma(\C)=\D, \quad \gamma(\D)=\A\B\C.$$

\noindent In the following we write $|w|$ for the length of a finite word $w$. Here are some useful properties of $\gamma$.

\begin{lemma}\label{lem:gamma} The morphism $\gamma$ has the following properties\\[.1cm]
{\rm i)\;} $|\gamma^n(\A)|=L_n$,  for all $n\ge 2$, where $L_n$ is the $n$th Lucas number (see next section).\\[.1cm]
{\rm ii)\;} $\gamma^n(\A)=\gamma^{n}(\C)$ and $\gamma^{n}(\A)=\gamma^{n+1}(\B)$ for all $n\ge 2$. 
\end{lemma}

\medskip

\proof {\rm i)\;} Starting at $n=2$, it follows easily with induction from the recursion of the Lucas numbers that one has
$|\gamma^n(\A)|=L_n,\; |\gamma^n(\B)|=L_{n-1},\; |\gamma^n(\C)|=L_n,\; |\gamma^n(\D)|=L_{n+1}$.\\[.1cm]
{\rm ii)\;} This follows immediately from $\gamma^2(\A)=\gamma(\A\B)=\A\B\C=\gamma(\D)=\gamma^2(\C)$.
 \hfill $\Box$

 \medskip


 \medskip

 It is notationally convenient to extend the semigroup of words to the free group of words. For example, one has $\D\C^{-1}\B^{-1}\B\C=\D$.

\section{Lucas numbers}\label{sec:luc}
The Lucas numbers $(L_n)=(2, 1, 3, 4, 7, 11, 18, 29, 47, 76,123, 199, 322,\dots)$ are defined by
$$  L_0 = 2,\quad L_1 = 1,\quad L_n = L_{n-1} + L_{n-2}\quad {\rm for \:}n\ge 2.$$
The Lucas numbers have a particularly simple base phi representation.

\noindent From  the well-known formula
$L_{2n}=\varphi^{2n}+\varphi^{-2n}$, and the recursion $L_{2n+1}=L_{2n}+L_{2n-1}$, we have for all $n\ge 1$
$$ \beta(L_{2n}) = 10^{2n}\cdot0^{2n-1}1,\quad \beta(L_{2n+1}) = 1(01)^n\cdot(01)^n.$$

\noindent {\bf Exercise} Show that the base phi representation of $L_{2n+1}+1$ equals $\beta(L_{2n+1}+1) = 10^{2n+1}\cdot(10)^n01$---see also Lemma 3.3. (2) in \cite{Hart99}, but note that these authors write the digits in reverse order.

\medskip

\noindent Since $\beta(L_{2n})$ consists of only 0's between the exterior 1's, the following lemma is obvious.

\begin{lemma}\label{lem:even}  For all $n\ge 1$ and $k=1,\dots,L_{2n-1}$
one has $ \beta(L_{2n}+k) =  \beta(L_{2n})+ \beta(k) = 10\dots0 \,\beta(k)\, 0\dots 01.$
\end{lemma}

 As in \cite{Hart98}, \cite{Hart99}, and \cite{San-San}, the strategy will be to partition the natural numbers in intervals $[L_n+1, L_{n+1}]$, and establish recursive relations for the $\beta$-expansions of the numbers in these intervals. However, an analogous formula as in Lemma \ref{lem:even} starting from an \emph{odd} Lucas number does not exist. To obtain recursive relations the interval $[L_{2n+1}+1, L_{2n+2}-1]$ has to be divided into three subintervals. These three intervals are
$$I_n:=[L_{2n+1}+1,\, L_{2n+1}+L_{2n-2}-1],\: J_n:=[L_{2n+1}+L_{2n-2},\, L_{2n+1}+L_{2n-1}], \: K_n:=[L_{2n+1}+L_{2n-1}+1,\, L_{2n+2}-1].$$
\noindent Note that $I_n$ and $K_n$ have the same length $L_{2n-2}-1$, that $J_n$ has length $L_{2n-3}+1$, and that the starting point $L_{2n+1}+L_{2n-2}$ of $J_n$  can be written as $2L_{2n}$.

From parts b.~and c.~of Proposition 3.1 and part c.~of Proposition 3.2 in the paper by Sanchis and Sanchis (\cite{San-San}) we obtain\footnote{N.B.: these authors write the beta-expansions in reverse order}  recursions for the beta-expansions of the natural numbers in the intervals $I_n$, $K_n$ and $J_n$.

\begin{lemma}\label{lem:odd}{\rm \bf(\cite{San-San})}  For all $n\ge 2$ and $k=1,\dots,L_{2n-2}-1$
$$ \beta(L_{2n+1}+k) = 1000(10)^{-1}\beta(L_{2n-1}+k)(01)^{-1}1001,$$ $$\beta(L_{2n+1}+L_{2n-1}+k)=1010(10)^{-1}\beta(L_{2n-1}+k)(01)^{-1}0001=10\beta(L_{2n-1}+k)(01)^{-1}0001.$$
Moreover, for all $n\ge 2$ and $k=0,\dots,L_{2n-3}$
$$\beta(L_{2n+1}+L_{2n-2}+k) = 10010(10)^{-1}\beta(L_{2n-2}+k)(01)^{-1}001001.$$
\end{lemma}

\medskip

As an illustration, we write out what Lemma~\ref{lem:odd} gives for $n=2$.
In the first part $k$ takes the values 1 and $L_2-1=2$, giving $(10)^{-1}\beta(5)(01)^{-1}=00\cdot 10$ and $(10)^{-1}\beta(6)(01)^{-1}=10\cdot 00$.
So the beta expansions of $L_5+1=12$, $L_5+2=13$, $L_5+L_3+1=16$ and $L_5+L_3+2=17$ are
$$\beta(12)=100000\cdot101001, \; \beta(13)=100010\cdot001001, \quad
\beta(16)=101000\cdot100001, \; \beta(17)=101010\cdot000001.$$
In the second part of Lemma~\ref{lem:odd} $k$ takes the values 0 and $L_1=1$, giving $(10)^{-1}\beta(3)(01)^{-1}=0\cdot$
and $(10)^{-1}\beta(4)(01)^{-1}=1\cdot$.
So the beta expansions of  $L_5+L_2+1=14$ and $L_5+L_2+1=15$ are
$$\beta(14)=100100\cdot001001, \; \beta(15)=100101\cdot001001.$$

\section{A proof of the conjectures}\label{sec:main}
	 The conjectures in the introduction will be part of the following more general result.

 \begin{theorem}  \label{th:d0d1} Let $\beta(N)=(d_i(N))$ be the base phi representation of a natural number $N$. Then:\\[.1cm]
  $d_0(N)=1$  \hspace*{1.1cm} if and only if $N=\lfloor n\varphi\rfloor +2n +1$ for some natural number $n$,\\[.1cm]
  $d_1d_0(N)=10$  \hspace*{.61cm}  if and only if $N=\lfloor n\varphi\rfloor +2n -1$ for some natural number $n$,\\[.1cm]
  $d_1d_0d_{-1}(N)=000$ if and only if $N=\lfloor n\varphi\rfloor +2n$ for some natural number $n$,\\[.1cm]
  $d_1d_0d_{-1}(N)=001$ if and only if $N=3\lfloor n\varphi\rfloor + n + 1$ for some natural number $n$.
\end{theorem}

\medskip

 It is convenient to code the four possibilities for the digits of $N$ by a map $T$ to an alphabet of four letters $\{\A,\B,\C,\D\}$.  We let\\[0.1cm]
 \hspace*{1.1cm} $T(N) =\A$ \;iff\;  $d_1d_0(N)=10$,\quad      $T(N)=\B$  \;iff\; $d_1d_0d_{-1}(N)=000$,  \\[0.1cm]
 \hspace*{1.1cm} $ T(N)=\C$ \;iff\; $d_0(N)=1$, \hspace*{0.65cm}  $T(N)=\D$   \;iff\;   $d_1d_0d_{-1}(N)=001$.\\[0.1cm]
 We thus have the following scheme.

 \bigskip

 \begin{tabular}{|r|c|c|}
   \hline
   $N^{\phantom{|}}$ & $\beta(N)$ & $T(N)$ \\[.0cm]
   \hline
   1 & \!\!\!$1$          & $\C^{\phantom{|}}$ \\
   2 & \:\,\,$10\cdot01$      & $\A$ \\
   3 & $100\cdot01$     & $\B$ \\
   4 & $101\cdot01$     & $\C$ \\
   5 & \:\,$1000\cdot1001$  & $\D$ \\
   6 & \:\,$1010\cdot0001$  & $\A$ \\
   7 & $10000\cdot0001$ & $\B$ \\
   8 & $10001\cdot 0001$  & $\C$\\
   \hline
 \end{tabular}  \quad
 \begin{tabular}{|r|c|c|}
   \hline
   $N^{\phantom{|}}$ & $\beta(N)$ & $T(N)$ \\[.0cm]
   \hline
   9  & \!\!\!$10010\cdot0101$     & $\A^{\phantom{|}}$ \\
   10 & \!\!\!$10100\cdot0101$     & $\B$ \\
   11 & \!\!\!$10101\cdot0101$     & $\C$ \\
   12 & $100000\cdot101001$  & $\D$ \\
   13 & $100010\cdot001001$  & $\A$ \\
   14 & $100100\cdot001001$  & $\B$ \\
   15 & $100101\cdot001001$  & $\C$ \\
   16 & $101000\cdot100001$  & $\D$ \\
   \hline
 \end{tabular} \quad
 \begin{tabular}{|r|c|c|}
   \hline
   $N^{\phantom{|}}$ & $\beta(N)$ & $T(N)$ \\[.0cm]
   \hline
   17  & \,\,\,$101010\cdot000001$    & $\A^{\phantom{|}}$ \\
   18 & $1000000\cdot000001$     & $\B$ \\
   19 & $1000001\cdot000001$     & $\C$ \\
   20 & $1000010\cdot010001$  & $\A$ \\
   21 & $1000100\cdot010001$  & $\B$ \\
   22 & $1000101\cdot010001$  & $\C$ \\
   23 & $1001000\cdot100101$  & $\D$ \\
   24 & $1001010\cdot000101$  & $\A$ \\
   \hline
 \end{tabular}

 \vspace*{0.8cm}

\noindent The reader may check the validity of the following $T$-values, which we use in the proof of Theorem \ref{th:rec}:
 $$ T(L_{2n})=\B,\;  T(L_{2n}+1)=\C,\; T(L_{2n+1}+1) =\D \quad{\rm for\:all}\: n\ge 1.\phantom{.............................}$$

\begin{theorem}  \label{th:ABCD}  The sequence $(T(N))_{N\ge 2}$ is the unique fixed point of the morphism $\gamma$.
\end{theorem}

\noindent Theorem \ref{th:ABCD} is an immediate consequence of Theorem \ref{th:rec}.

\begin{theorem}  \label{th:rec} Let $\gamma$ be the morphism given by $\A\mapsto\A\B,\, \B\mapsto\C,\,\C\mapsto\D,\,\D\mapsto\A\B\C.$  Then\\[.1cm]
{\rm  a)} $T(2)\,T(3)\cdots T(L_n\!+\!1) = \gamma^n(\A)$ for $n\ge 2$\\[.1cm]
{\rm  b)} $T(L_n\!+\!2)\,T(L_n\!+\!3)\cdots T(L_{n+1}\!+\!1) = \gamma^{n-1}(\A)$ for $n\ge 3$.
\end{theorem}

\proof We prove a) and b) simultaneously by induction.\\[.1cm]
 For $n=2$, $L_2=3$, and
 one finds $T(2)T(3)T(4)=\A\B\C$, which indeed equals $\gamma^2(\A)$.\\[.1cm]
  Also for $n=3$, one has   $T(2)T(3)T(4)T(5)=\A\B\C\D=\gamma^3(\A)$.\\[.1cm]
  Part b) for $n=3$ is checked by $T(6)T(7)T(8)=\A\B\C=\gamma^2(\A)$.\\[.1cm]
In the following we do not formally perform an induction step $n\rightarrow n+1$, but show how $T$-images of intervals can be expressed in $T$-images of intervals with lower indices. We have for part a)
\begin{align*}
 T(2)\cdots T(L_{n+1}\!+\!1)& =  T(2)\cdots T(L_n\!+\!1)\,T(L_n\!+\!2)\cdots T(L_{n+1}\!+\!1)\\
                             & =  \gamma^n(\A)\,\gamma^{n-1}(\A)\\
                             & =  \gamma^{n}(\A\B) =\gamma^{n+1}(\A).
\end{align*}
Here we used Lemma \ref{lem:gamma} part ii).\\[.1cm]
For part b), this formula follows for  even indices directly from  Lemma \ref{lem:even} and part a):
\begin{align*}
T(L_{2n}\!+\!2)\cdots T(L_{2n+1})\, T(L_{2n+1}\!+\!1)& =T(L_{2n}\!+\!2)\cdots T(L_{2n+1})\, \D \\
     & =T(2)\dots T(L_{2n-1})\,\D\\ & =T(2)\dots T(L_{2n-1})\,T(L_{2n-1}+1)= \gamma^{2n-1}(\A).
\end{align*}
For  odd indices, we use Lemma \ref{lem:odd}. We have
\begin{align*}
 T(L_{2n+1}+1) \cdots T(L_{2n+1}+L_{2n-2}\!-\!1)  &= T(L_{2n+1}+1)\,\gamma^{2n-2}(\A)\,T(L_{2n}+1)^{-1}\,T(L_{2n})^{-1} \\
                                  & =\D\,\gamma^{2n-2}(\A)\,\C^{-1}\B^{-1},\\
T(L_{2n+1}+L_{2n-2})\cdots T(L_{2n+1}+L_{2n-1})   &=  T(L_{2n-2})\,T(L_{2n-2}+1)\,\cdots T(L_{2n-1}+1)\,T(L_{2n-1}+1)^{-1} \\
                                  & =\B\,\C\, \gamma^{2n-3}(\A)\,\D^{-1},\\
T(L_{2n+1}+L_{2n-1}+1) \cdots T(L_{2n+2}\!-\!1) &= \D\,\gamma^{2n-2}(\A)\,\C^{-1}\B^{-1}.\\
\phantom{noline}
\end{align*}

 \vspace*{-.5cm}

\noindent Concatenating the $T$-images of the intervals $I_n, J_n$ and $K_n$, we obtain, using Lemma \ref{lem:gamma} part ii)
\begin{align*}
& T(L_{2n+1}+2) \cdots T(L_{2n+2}+1)  =\\
& T(L_{2n-1}+1)^{-1}\,\D\,\gamma^{2n-2}(\A)\,\C^{-1}\B^{-1}\B\,\C\, \gamma^{2n-3}(\A)\,\D^{-1}\D\,\gamma^{2n-2}(\A)\,\C^{-1}\B^{-1}\B\C=\\
& \gamma^{2n-2}(\A)\, \gamma^{2n-3}(\A)\,\gamma^{2n-2}(\A)\,=\, \gamma^{2n-2}(\A\B\C) = \gamma^{2n-2}(\gamma^2(\A))\,=\,\gamma^{2n}(\A).\hspace*{7.5cm} \Box
\end{align*}\vspace*{-.1cm}

 \medskip

 \noindent {\it Proof of Theorem \ref{th:d0d1}:}
 From Theorem \ref{th:ABCD} we know that the digit $d_0(N)=1$ iff $T(N)=\C$, where (with some abuse of notation) $T=\C\A\B\C\A\B\C\D\dots$ is the fixed point of $\gamma$, prefixed by $\C$.  We see from the form of $\gamma^2$ that (apart from the prefix $\C$) $T$ is a concatenation of the words $\A\B\C$   and $\D$.  Suppose we apply a code:  $\psi(\A\B\C)=0$,  $\psi(\D)=1$.
  Then $\gamma$ induces a morphism $\sigma$ on the alphabet $\{0,1\}$:
  $$\sigma:\quad 0\mapsto \psi(\gamma(\A\B\C))=\psi(\A\B\C\D)=01, \quad 1\mapsto \psi(\gamma(\D)) =\psi(\A\B\C)=1.$$
  We see that $\sigma$ is the Fibonacci morphism, with fixed point $x_{0,1}$.
  But the 0's in $x_{0,1}$ occur at positions $\lfloor n\varphi\rfloor$, $n=1,2\dots$ (see, e.g., \cite{lothaire}). Since the differences between the indices of the positions of $\C$ in $T$ are expanded by 2 by the inverse of $\psi$, and because of the prefix $\C$, this implies that the $\C$'s occur at positions $\lfloor n\varphi\rfloor+2n+1$, for $n=0,1,\dots$. But obviously $\A$'s always occur at two places before a $\C$, implying that the positions of $\A$ are given by $\lfloor n\varphi\rfloor+2n-1$, for $n=1,\dots$.
  Similarly the positions of $\B$ are given by $\lfloor n\varphi\rfloor+2n$.

  Finding the  positions of $\D$ is more involved. Consider the locations of $\D$ in the morphism $\gamma^4$:
  $$\gamma^4: \quad \A\mapsto \A\B\C\underline{\D}\A\B\C,\; \B\mapsto \A\B\C\underline{\D}, \;
  \C\mapsto \A\B\C\underline{\D}\A\B\C, \; \D\mapsto \A\B\C\underline{\D}\A\B\C\A\B\C\underline{\D}.$$
 We see from this that the difference between the indices of occurrence of $\D$ in $T=\gamma^4(T)$ is always 4 or 7.
 Moreover, the distances generated by $\A,\B,\C$ and $\D$ under $\gamma$ are respectively 7, 4, 7, and the pair 7,4.
 Mapping $\A\mapsto 7, \B \mapsto 4, \C \mapsto 7, \D \mapsto 74$, the morphism $\gamma$ induces for $\A,\C$ and $\B$ a morphism $7\mapsto 74, 4\mapsto 7$. Moreover, this morphism is compatible with the part induced by $\D$: \: $74\mapsto 747$. It follows that the sequence of differences of  indices of occurrence of $\D$ is nothing else but the Fibonacci sequence $x_{7,4}$ on the alphabet $\{7,4\}$. Lemma \ref{lem:diff} then gives that this sequence equals $(3\lfloor n\varphi\rfloor + n + 1)_{n\ge 1}$. \hfill\endpf

 \smallskip

 \begin{remark} With induction, using Lemma \ref{lem:even} and \ref{lem:odd},  one proves that $d_1d_0(N)=10$ forces $d_{-1}(N)=0$. It follows that Theorem \ref{th:d0d1} implies that

\qquad Digit $d_{-1}(N)=1$ if and only if $N=3\lfloor n\varphi\rfloor + n + 1$ for some natural number $n$.
\end{remark}


\section{A general result}\label{sec:general}
	
Here we given an expression for the set of $N$ with $d_k(N)=1$ for any $k>1$. Recall that we partitioned the natural numbers in Lucas intervals\: $\Lambda_{2n}:=[L_{2n},L_{2n+1}]$ and $\Lambda_{2n+1}:=[L_{2n+1}+1, L_{2n+2}-1]$.\\ The basic idea behind this partition is that if
 $$\beta(N) = d_{L}d_{L-1}\dots d_1d_0\cdot d_{-1}d_{-2} \dots d_{R-1}d_R,$$
then the left most index $L=L(N)$ and the right most index $R=R(N)$ satisfy
$$L(N)=2n\!+1, \,R(N)=2n \;{\rm iff}\; N\in \Lambda_{2n}, \quad L(N)=2n\!+2= R(N) \;{\rm iff}\; N\in \Lambda_{2n+1}.$$
This is not hard to see from the simple expressions we have for the $\beta$-expansions of the Lucas numbers, see also Theorem 1 in \cite{Grabner94}. For the  cardinality $|\Lambda_{n}|$ of  $\Lambda_{n}$ we have (of course!)
$$|\Lambda_{n}| = \lfloor \varphi^{n+1} \rfloor-\lfloor \varphi^n \rfloor.$$
Note that we also have $|\Lambda_{2n}|=L_{2n-1}+1$, and $|\Lambda_{2n+1}|=L_{2n}-1$, the expressions used in \cite{San-San}.
It can therefore be checked easily that our Theorem \ref{th:dN} implies the main result of \cite{San-San} (for positive $k$).

 \begin{theorem}  \label{th:dN} Let $\beta(N)=(d_i(N))$ be the base phi representation of a natural number $N$, and let $k\ge 2$. Then
  $d_k(N)=1$  if and only if $N$ is a member of one of the
 generalized Beatty sequences $(\lfloor n\varphi\rfloor L_k +nL_{k-1} +r)$, where  $r=r_1,\,r_1\!+1,\dots,r_1\!+|\Lambda_{k}|\!-1$, with $r_1=-L_{k-1}$ if $k$ is even, and $r_1=-L_{k-1}\!+\!1$ if $k$ is odd.
\end{theorem}

\proof It turns out that the coding with the alphabet $\{\A,\B,\C,\D\}$ is still useful. In fact, we extend this alphabet to an alphabet $\{\A_0,\A_1,\B_0,\B_1,\C_0,\C_1,\D_0,\D_1\}$ via the extended coding $T_+$ defined for $j=0,1$ by
$$T_+(N)= \A_j {\;\;\rm iff\;\;}  d_k(N)=j,\: T(N)=\A,\; \dots,\; T_+(N)= \D_j {\;\;\rm iff\;\;}  d_k(N)=j,\: T(N)=\D.   $$
We also want to extend the morphism $\gamma$ to a morphism $\gamma_+$. Here it turns out that one has to extend $\gamma^{k+2}$ instead of $\gamma$. For simplicity in notation we suppress the dependence on $k$ in $\gamma_+$. We obtain $\gamma_+$ by looking at  $\gamma^{k+2}(\A)\gamma^{k+2}(\B)\gamma^{k+2}(\C)\gamma^{k+2}(\D)$---note that this word is always a prefix of $(T(N))_{N\ge 2}$ as a consequence of Theorem \ref{th:ABCD}. We define
\begin{align*}
\gamma_+(\A_0)&= \gamma_+(\A_1) = T_+(2)\dots T_+(L_{k+2}+1),\\
\gamma_+(\B_0)&= \gamma_+(\B_1) = T_+(L_{k+2}+2)\dots T_+(L_{k+2}+L_{k+1}+1)= T_+(L_{k+2}+2)\dots T_+(L_{k+3}+1),\\
\gamma_+(\C_0)&= \gamma_+(\C_1) = T_+(L_{k+3}+2)\dots T_+(L_{k+3}+L_{k+2}+1)= T_+(L_{k+3}+2)\dots T_+(L_{k+4}+1),\\
\gamma_+(\D_0)&= \gamma_+(\D_1) = T_+(L_{k+4}+2)\dots T_+(L_{k+4}+L_{k+3}+1)=T_+(L_{k+4}+2)\dots T_+(L_{k+5}+1).
\end{align*}

\noindent In view of the complexity of the proof we start with the case $k=2$, so $\gamma^{k+2}=\gamma^4$, and $\gamma_+$ has the form:
\begin{align*}
\gamma_+(\A_0)&= \gamma_+(\A_1) = \A_0\B_1\C_1\D_0\A_0\B_0\C_0,\\
\gamma_+(\B_0)&= \gamma_+(\B_1) = \A_0\B_1\C_1\D_0,\\
\gamma_+(\C_0)&= \gamma_+(\C_1) = \A_0\B_1\C_1\D_0\A_0\B_0\C_0,\\
\gamma_+(\D_0)&= \gamma_+(\D_1) = \A_0\B_1\C_1\D_0\A_0\B_0\C_0\A_0\B_1\C_1\D_0.
\end{align*}

\noindent Here the $\B_1\C_1$ in $\gamma_+(\A_j)$ is coming from the first couple of 1's in $d_2(N)$ occurring in $\Lambda_2 = [L_2,L_3] = [3,4].$

\medskip

We claim that   $(T_+(N))_{N\ge 2}$ is the unique fixed point of $\gamma_+$.
We will prove this in a way similar to the proof of Theorem \ref{th:rec}.

\medskip

\noindent CLAIM:\\
$\boxplus$ {\rm  a)} $T_+(2)\cdots T_+(L_{4n}\!+\!1) = \gamma_+^n(\A_0)$ for $n\ge 1$\\[.1cm]
$\boxplus$ {\rm  b)} $T_+(L_{4n}\!+\!2)\cdots T_+(L_{4n+1}\!+\!1) = \gamma_+^{n}(\B_0)$ for $n\ge 1$.\\[.1cm]
$\boxplus$ {\rm  c)} $T_+(L_{4n+1}\!+\!2)\cdots T_+(L_{4n+2}\!+\!1) = \gamma_+^{n}(\C_0)$ for $n\ge 1$.\\[.1cm]
$\boxplus$ {\rm d)} $T_+(L_{4n+2}\!+\!2)\cdots T_+(L_{4n+3}\!+\!1) = \gamma_+^{n}(\D_0)$ for $n\ge 1$.\\[.1cm]
$\boxplus$ {\rm  e)} $T_+(L_{4n+3}\!+\!2)\cdots T_+(L_{4n+4}\!+\!1) = \gamma_+^{n}(\A_0\B_0\C_0)$ for $n\ge 1$.

\medskip

\noindent {\it Proof of the claim:} This will be done with induction, with an unexpected twist.\\[.1cm]
First the case $n=1$. \\[.1cm]   
By definition one has $\boxplus$ {\rm  a)} $T_+(2)\cdots T_+(L_{4}\!+\!1) = \gamma_+(\A_0)$,
$\boxplus$ {\rm  b)} $T_+(L_{4}\!+\!2)\cdots T_+(L_5\!+\!1) = \gamma_+(\B_0)$,
$\boxplus$ {\rm  c)} $T_+(L_{5}\!+\!2)\cdots T_+(L_6\!+\!1) = \gamma_+(\C_0)$, and
$\boxplus$ {\rm  d)} $T_+(L_{6}\!+\!2)\cdots T_+(L_{7}\!+\!1) = \gamma_+(\D_0)$.  \\[.1cm]
What remains is
$\boxplus$ {\rm  e)} $T_+(L_{7}\!+\!2)\cdots T_+(L_{8}\!+\!1) = \gamma_+(\A_0\B_0\C_0)$, which can be proved by using Lemma \ref{lem:odd}:\\[.1cm]
the central part of $\beta(L_{7}\!+\!k)$ equals $\beta(L_{5}\!+\!k)$ for $k=1,\dots L_4-1$,
yielding $T_+(L_{7}\!+\!2)\cdots T_+(L_{7}\!+\!L_4\!-\!1) =  \gamma_+(\C_0)\C_0^{-1}\B_0^{-1}$. Similarly,\:
$T_+(L_{7}\!+\!L_5\!+\!1)\cdots T_+(L_{8}\!-\!1) = \D_0 \gamma_+(\C_0)\C_0^{-1}\B_0^{-1}$. In between we have
$T_+(L_{7}\!+\!L_4)\cdots T_+(L_{7}\!+\!L_4\!+\!L_3) = \B_0\C_0 \gamma_+(\B_0)\D_0^{-1}$. Pasting these three words  together, and adding the two letters $T_+(L_{8}) = \B_0$, and   $T_+(L_{8}\!+\!1) = \C_0$, we obtain the word $\gamma_+(\C_0\B_0\C_0)=\gamma_+(\A_0\B_0\C_0)$.

\medskip

\noindent Next we make the induction step $n\rightarrow n+1$.

\medskip


\noindent$\boxplus$ {\rm  a)} Here one splits $T_+(2)\cdots T_+(L_{4(n+1)}\!+\!1)$ into 5
subwords $T_+(L_{4n+j}+2)\cdots T_+(L_{4n+j+1}+1)$, $j=0,\dots,4$. The induction hypothesis then  gives
$$T_+(2)\cdots T_+(L_{4(n+1)}\!+\!1)=
\gamma_+^n(\A_0)\gamma_+^n(\B_0)\gamma_+^n(\C_0)\gamma_+^n(\D_0)\gamma_+^n(\A_0\B_0\C_0)=\gamma_+^{n+1}(\A_0).$$
\noindent $\boxplus$ {\rm  b)} From Lemma \ref{lem:even} one obtains from the induction hypothesis, again with a splitting
$$T_+(L_{4(n+1)}\!+\!2)\cdots T_+(L_{4(n+1)+1}\!+\!1)=T_+(2)\cdots T_+(L_{4n+3}\!+\!1)=
       \gamma_+^n(\A_0)\gamma_+^n(\B_0)\gamma_+^n(\C_0)\gamma_+^n(\D_0)=\gamma_+^{n+1}(\B_0).$$
\noindent$\boxplus$  {\rm  c)} This is more involved, as we have to use  Lemma \ref{lem:odd}. This lemma yields
\begin{align*}
T_+(L_{4(n+1)+1}\!+\!2)\cdots T_+(L_{4(n+1)+1}\!+\!L_{4n+2}-1)&=T_+(L_{4(n+1)-1}\!+\!2)\cdots T_+(L_{4(n+1)-1}\!+\!L_{4n+2}-1)\\
   &=T_+(L_{4n+3}\!+\!2)\cdots T_+(L_{4n+4}-1) = \gamma_+^n(\A_0\B_0\C_0)\C_0^{-1}\B_0^{-1},
\end{align*}
where we used part e) of the induction hypothesis in the last step. For the `middle part' Lemma \ref{lem:odd} yields
\begin{align*}
 T_+(L_{4(n+1)+1}\!+\!L_{4n+2})\cdots T_+(L_{4(n+1)+1}\!+\!L_{4n+3})&=T_+(L_{4n+2})\cdots T_+(L_{4n+3})=\B_0\C_0\gamma_+^n(\D_0)\D_0^{-1}
\end{align*}
The last part is similar to the first part. Pasting the three parts together, and adding $\B_0\C_0$ at the end we obtain
\begin{align*}
T_+(L_{4(n+1)+1}\!+\!2)\cdots T_+(L_{4(n+1)+2}\!+\!1)&=
 \gamma_+^n(\A_0\B_0\C_0)\C_0^{-1}\B_0^{-1}\B_0\C_0\gamma_+^n(\D_0)\D_0^{-1}\D_0\gamma_+^n(\A_0\B_0\C_0)\C_0^{-1}\B_0^{-1}\B_0\C_0\\
 &=\gamma_+^n(\A_0\B_1\C_1)\gamma_+^n(\D_0)\gamma_+^n(\A_0\B_0\C_0)=\gamma_+^{n+1}(\C_0).
\end{align*}
\noindent$\boxplus$  {\rm  d)} From Lemma \ref{lem:even} one obtains
\begin{align*}
T_+(L_{4(n+1)+2}\!+\!2)\cdots T_+(L_{4(n+1)+3}\!+\!1)&=T_+(2)\cdots T_+(L_{4n+5}\!+\!1)\\
&= T_+(2)\cdots T_+(L_{4n+4}\!+\!1)T_+(L_{4n+4}\!+\!2)\cdots T_+(L_{4n+5}\!+\!1)\\   &=\gamma_+^{n+1}(\A_0)\,\gamma_+^{n+1}(\B_0)=\gamma_+^{n+1}(\D_0).
\end{align*}
Here we could not use the induction hypothesis, but can apply part a) and b) already proved above.\\[0.1cm]
\noindent$\boxplus$  {\rm  e)} Again, we have to use  Lemma \ref{lem:odd}. This lemma yields
\begin{align*}
T_+(L_{4(n+1)+3}\!+\!2)\cdots T_+(L_{4(n+1)+3}\!+\!L_{4n+2}-1)&=T_+(L_{4(n+1)+1}\!+\!2)\cdots T_+(L_{4(n+1)+1}\!+\!L_{4n+4}-1)\\
   &=T_+(L_{4n+5}\!+\!2)\cdots T_+(L_{4n+6}-1) = \gamma_+^{n+1}(\C_0)\C_0^{-1}\B_0^{-1},
\end{align*}
where we used part c) already proved above. For the 'middle part' Lemma \ref{lem:odd} yields
\begin{align*}
 T_+(L_{4(n+1)+3}\!+\!L_{4n+4})\cdots T_+(L_{4(n+1)+3}\!+\!L_{4n+5})&=T_+(L_{4n+4})\cdots T_+(L_{4n+5})=\B_0\C_0\gamma_+^{n+1}(\B_0)\D_0^{-1},
\end{align*}
where we used part b) already proved above. \\ The last part is similar to the first part. Pasting the three parts together we obtain
$$ T_+(L_{4(n+1)+3}\!+\!2)\cdots T_+(L_{4(n+1)+4}\!+\!1)= \gamma_+^{n+1}(\C_0)\gamma_+^{n+1}(\B_0)\gamma_+^{n+1}(\C_0)
 =\gamma_+^{n+1}(\A_0\B_0\C_0).$$
This finishes the proof of the claim. To finish the proof of the theorem for the case $k=2$, we note that the situation is almost
identical\footnote{This observation also leads to a more or less independent proof of Theorem \ref{th:dN} for $k=2$: $\B_1\C_1$ occurs always immediately before $\D_0$, so the positions of $\B_1$, respectively $\C_1$, are just those of $\D$ in Theorem \ref{th:d0d1} shifted by -1 and -2.} to the appearance of $\D$ in $\gamma^4(\A),\dots,\gamma^4(\D)$ at the end of the proof of Theorem \ref{th:ABCD}: the words $\B_1\C_1$ occur at indices which differ by $7$ or $4$, and these differences occur as $x_{7,4}$, the Fibonacci word on the alphabet $\{7,4\}$. An application of Lemma \ref{lem:diff} then gives that the numbers $N$ with $d_2(N)=1$ occur as $N=3\lfloor n\varphi\rfloor  +n +r$ with two possibilities for  $r$, which are found to be $r=0$ and $r=-1$.

\bigskip

Consider in general the case of an even integer $2k,\, k=1,2\dots.$ One first proves that  $(T_+(N))_{N\ge 2}$ is the unique fixed point of $\gamma_+$, following the same scheme as in the proof for the $k=2$ case. Next, one has to sort out where the $N$ with $d_{2k}(N)=1$ appear with respect to the $\gamma_+(\A_0),\dots,\gamma_+(\D_0)$ in the fixed point of $\gamma_+$. The first time $d_{2k}(N)=1$ appears is for $N=L_{2k}$, the first number  in $\Lambda_{2k}$, and all other $N$ in  $\Lambda_{2k}$ also have $d_{2k}(N)=1$. By Lemma \ref{lem:even}, these trains of $N$'s with $d_{2k}(N)=1$ also appear at the end of $\Lambda_{2k+2}$ (excepting $N=L_{2k+3}+1$). Since they can not appear in $\Lambda_{2k+1}$, this \emph{is} the second appearance of the train. Application of Lemma \ref{lem:odd}, and another time Lemma \ref{lem:even}, then gives that the third appearance is in $\Lambda_{2k+3}$, and the fourth and fifth appearance are in $\Lambda_{2k+4}$.
Moreover, these three Lucas intervals correspond---except for one or two symbols at the begin and at the end---to the intervals used to define $\gamma_+(\B_0)$, $\gamma_+(\C_0)$, and $\gamma_+(\D_0)$, and at the same time it shows that $\gamma_+(\C_0)=\gamma_+(\A_0)$, and $\gamma_+(\D_0)=\gamma_+(\B_0)\gamma_+(\C_0)$.

 This means that the situation is very much like the appearance of $\B_1\C_1$ in the words $\gamma_+(\A_0),\dots,\gamma_+(\D_0)$ in the $k=2$ case treated above: the trains occur at indices which differ by $L_{2k+2}$ or $L_{2k+1}$, and these differences occur as $x_{L_{2k+2}, L_{2k+1}}$, the Fibonacci word on the alphabet $\{L_{2k+2}, L_{2k+1}\}$. An application of Lemma \ref{lem:diff} then gives that the numbers $N$ in the train occur as $\lfloor n\varphi\rfloor L_{2k} +nL_{2k-1} +r$ for some $r$, since
$$L_{2k+2}- L_{2k+1}=L_{2k},\quad {\rm and \quad} 2L_{2k+1}-L_{2k+2}=L_{2k-1}.$$
Substituting $n=1$, corresponding to the first train, with first element $N=L_{2k}$, gives $r_1=-L_{2k-1}$. The length of the train is of course $|\Lambda_{2k}|$.

\medskip

The proof for odd integers $2k+1$ follows the same steps, the sole difference being that $r_1$ turns out to be one larger, due to the fact that $\Lambda_{2k}$ starts at $L_{2k}$, but $\Lambda_{2k\!+\!1}$ starts at $L_{2k\!+\!1}\!+\!1$ . \endpf

\medskip

 \begin{remark} A result similar to Theorem \ref{th:dN} will hold for digits $d_N(k)$ with $k$ negative, but the situation is somewhat more complex. One has, for example,

\; Digit $d_{-2}(N)=1$ if and only if $N=4\lfloor n\varphi\rfloor +3n + r$ for some $r=2,3,4$ and some non-negative integer $n$.
\end{remark}

\end{document}